\documentclass[11pt,twoside]{article}

\usepackage{a4wide}
\usepackage{amsfonts}
\usepackage{amssymb}
\usepackage{amsmath}
\usepackage{graphicx}

\newcommand{\ignore}[1]{}

\def\@begintheorem#1#2{\par\bgroup{\sc #1\ #2. }\it\ignorespaces}
\def\@opargbegintheorem#1#2#3{\par\bgroup{\sc #1\ #2\ (#3). } \it\ignorespaces}
\def\@endtheorem{\egroup}
\newtheorem{theorem}{Theorem}[section]
\newtheorem{corollary}[theorem]{Corollary}
\newtheorem{lemma}[theorem]{Lemma}

\newtheorem{example}[theorem]{Example}
\newtheorem{proposition}[theorem]{Proposition}
\newtheorem{definition}[theorem]{Definition}
\newcommand{\bt}[1]{\begin{theorem}\label{#1}}
\newcommand{\bc}[1]{\begin{corollary}\label{#1}}
\newcommand{\bl}[1]{\begin{lemma}\label{#1}}
\newcommand{\be}[1]{\begin{example}\label{#1}}
\newcommand{\bp}[1]{\begin{proposition}\label{#1}}
\newcommand{\ba}[1]{\begin{algorithm}\rm\label{#1}}
\newcommand{\bd}[1]{\begin{definition}\rm\label{#1}}{\normalsize }
\newcommand{\bpr}{\noindent {\em Proof. }}
\newcommand{\et}{\end{theorem}}
\newcommand{\ec}{\end{corollary}}
\newcommand{\el}{\end{lemma}}
\newcommand{\ee}{\end{example}}
\newcommand{\ep}{\end{proposition}}
\newcommand{\ed}{\end{definition}}
\newcommand{\epr}{{\ \vbox{\hrule\hbox{%
\vrule height1.3ex\hskip0.8ex\vrule}\hrule}}\\\par}
\newcommand{\mepr}{{\ \ \ \vbox{\hrule\hbox{%
\vrule height1.3ex\hskip0.8ex\vrule}\hrule}}}

\def\R{\mathbb{R}}
\def\Z{\mathbb{Z}}

\def\sh{\overline}
\def\cc{{\overline c}}
\def\xx{{\overline x}}

\begin{document}

\title{\bf Approximate Shifted Combinatorial Optimization}

\author{
Martin Kouteck\'y
\thanks{\small Charles University, Prague.
Email: koutecky@kam.mff.cuni.cz}
\and
Asaf Levin
\thanks{\small Technion - Israel Institute of Technology, Haifa.
Email: levinas@ie.technion.ac.il}
\and
Syed M. Meesum
\thanks{\small Institute of Mathematical Sciences, HBNI, Chennai.
Email: meesum@imsc.res.in}
\and
Shmuel Onn\footnotemark[2]
\thanks{\small Technion - Israel Institute of Technology, Haifa.
Email: onn@ie.technion.ac.il}
}
\date{}

\maketitle

\begin{abstract}
{\em Shifted combinatorial optimization} is a new nonlinear optimization
framework, which is a broad extension of standard combinatorial optimization,
involving the choice of several feasible solutions at a time. It captures well studied and diverse problems ranging from congestive to partitioning problems.
In particular, every standard combinatorial optimization problem has its shifted counterpart, which is typically much harder. Here we initiate a study of
approximation algorithms for this broad optimization framework.

\vskip.2cm
\noindent {\bf Keywords:} combinatorial optimization, packing problems,
independence system, matroid, matching, approximation algorithms, linear optimization.

\end{abstract}

\section{Introduction}

The following optimization problem has been studied extensively in the literature.
\vskip.2cm\noindent{\bf (Standard) Combinatorial Optimization.}
Given $S\subseteq\{0,1\}^d$ and $w\in\Z^d$, solve
\begin{equation}\label{standard}
\max\{ws\ :\ s\in S\}\ .
\end{equation}
The complexity of the problem depends on $w$ and the type and presentation of $S$.
Often, $S$ arises as the set of indicating vectors of members of a family of subsets over a ground set $[d]:=\{1,\dots,d\}$ such as the family of matchings in a given graph with $d$ edges
or the set of independent sets in a matroid over $[d]$ given by an independence oracle.
See \cite{Sch} for a detailed account of the literature and bibliography of thousands
of articles on this problem.

In this article we study a broad nonlinear extension of this problem,
which involves the choice of several feasible solutions from $S$ at a time,
and which is defined as follows.

We denote the $i$-th row and $j$-th column of a matrix $x$ by $x_i$ and $x^j$ respectively.
For a set $S\subseteq\R^d$ let $S^n$ be the set of $d\times n$
matrices having each column in $S$,
$$S^n\ :=\ \{x\in\R^{d\times n}\ :\ x^j\in S\,,\ j=1,\dots,n \}\ .$$
Call matrices $x,y\in\R^{d\times n}$ {\em equivalent} and write $x\sim y$ if each row of $x$
is a permutation of the corresponding row of $y$. The {\em shift} of $x\in\R^{d\times n}$
is the unique matrix $\xx\in\R^{d\times n}$ satisfying $\xx\sim x$ and
$\xx^1\geq \cdots\geq \xx^n$, that is, the unique matrix equivalent to $x$ with each
row nonincreasing.  We say that a matrix $x$ is {\em shifted} if $x= \xx$.
We study the following broad nonlinear combinatorial optimization problem
(also called the {\em shifting problem} and denoted as SCO).

\vskip.2cm\noindent{\bf Shifted Combinatorial Optimization (SCO).}
Given $S\subseteq\{0,1\}^d$ and $c\in\Z^{d\times n}$, solve
\begin{equation}\label{shift}
\max\{c\xx\ :\ x\in S^n\}\ .
\end{equation}

\vskip.2cm
This problem has a very broad expressive power. In particular, every standard
combinatorial optimization has its shifted counterpart which is typically much harder.
For instance, when $S$ is given explicitly as a list of vectors, the standard problem is
trivial, but the shifted counterpart may be hard. To see this, let $G$ be a graph with
$d$ edges and $m$ vertices. Let $S:=\{s^1,\dots,s^m\}\subseteq\{0,1\}^d$ with $s^j$ the
indicator of the set of edges incident on vertex $j$. Define $c\in\Z^{d\times n}$ by
$c_{i,1}:=0$ for all $i$ and $c_{i,j}:=-1$ for all $i$ and $j\geq 2$. Then the optimal
objective value of the shifted problem is $0$ if and only if $G$ has an independent
set of size $n$.

One interpretation of the problem is in terms of the social cost in a {\em congestion game}
\cite{Ros}. We are given a set $S\subseteq\{0,1\}^d$ of indicators of members of a family
over $[d]$ as above. For $i=1,\dots,d$ we are now given a function
$f_i:\{0,1,\dots,n\}\rightarrow \Z$. Each $x\in S^n$ represents a choice of $n$ players with
$x^j\in S$ the choice of player $j$. The {\em congestion} of $x$ at element $i\in[d]$ is
the number $\sum_{j=1}^n x_{i,j}$ of players using $i$ in $x$. The cost of $x$ at $i$ is
the value $f_i(\sum_{j=1}^n x_{i,j})$ of $f_i$ on the congestion at $i$. The social cost
of $x$ is $\sum_{i=1}^d f_i(\sum_{j=1}^n x_{i,j})$ and we want to find $x\in S^n$
minimizing social cost. For instance, $S$ may be the set of $s-t$ dipaths in a digraph;
each player chooses a dipath; and the cost at edge $i$ may be an increasing function $f_i$
of the congestion at $i$. Now define $c\in\Z^{d\times n}$ by $c_{i,j}:=f_i(j-1)-f_i(j)$
for all $i$ and $j$. Then for every $x\in S^n$ we have
$c\xx=\sum_{i=1}^df_i(0)-\sum_{i=1}^d f_i(\sum_{j=1}^n x_{i,j})$ and so $x$ maximizes
$c\xx$ if and only if it minimizes the social cost. So these congestion and shifted
problems are equivalent. In particular, a case of special interest is when the $f_i$
are {\em convex}, implying that the rows of $c$ are nonincreasing, that is, $c=\cc$
is {\em shifted}. This case will be considered in Section \ref{constant_appx}.

Shifted combinatorial optimization has been introduced and studied recently in
\cite{GHKO,KOS,LO}. In \cite{KOS} it was shown that for $S=\{s\in\{0,1\}^d\,:\,As\leq b\}$
presented by linear inequalities where $A$ is a totally unimodular matrix and $b$ is
integer, with shifted $c$, the shifting problem can be solved in polynomial time.
This in particular implies that shifting with shifted $c$ can be solved in polynomial time
for $S$ the set of matchings in a bipartite graph or the set of $s-t$ dipaths in a digraph.
In \cite{LO} it was shown that with shifted $c$, for $S$ the set of independent sets in a
matroid presented by a linear optimization oracle or $S$ the intersection of two so-called
strongly-base-orderable matroids, shifting with shifted $c$ can again be done in polynomial time.
Finally, in \cite{GHKO} the parameterized complexity of the problem with $S$ given explicitly
and $|S|$ as parameter was studied, and it was also shown that shifted combinatorial
optimization parameterized by $\mu$ and $\tau$ over any property
defined by any monadic-second-order-logic formula of length at most $\mu$
over any graph of tree-width at most $\tau$, is in the complexity class XP.

In this article we initiate a study of approximation algorithms for shifted
combinatorial optimization. We assume that $S$ is an {\em independence system}, also
called {\em downward monotone}, namely, $0\in S$ and if $u\in\{0,1\}^d$ and
$u\leq v\in S$ then $u\in S$. We assume that $S$ is presented by a
{\em linear optimization oracle} that, queried on $w\in\Z^d$, returns $s\in S$
maximizing $ws$, that is, we can at least solve the standard
combinatorial optimization problem over $S$ to begin with.
For instance, $S$ may be the independent sets in a matroid or the matchings in a graph.

\vskip.2cm
The article is organized as follows. In Section \ref{matchings} we discuss the complexity of
the problem for matchings. This is done via the so-called {\em prescribed congestion problem}
which is of interest on its own right. We show that for shifted $c$, shifting
is hard already for $n=2$ and matchings in general graphs, and for general $c$,
it is hard already for $n=2$ and matchings in bipartite graphs.
This motivates the need in approximation algorithms for shifted combinatorial optimization
over independence systems, which are developed in Sections \ref{constant_appx} -- \ref{small}.
As usual, an algorithm for a maximization problem has {\em approximation ratio $\alpha$}
if it always returns a feasible solution with objective value that is
at least $\alpha$ times the optimal objective value (that is, we use the convention
of $\alpha \leq 1$ for maximization problems). Also, an algorithm is {\em polynomial time}
if its running time including the number of calls to the oracle presenting $S$ is polynomial.

\vskip.2cm
In Section \ref{constant_appx} we prove the following
result (see therein for the precise statement).

\vskip.2cm\noindent{\bf Theorem \ref{constant_ratio}.}
SCO for any shifted $c$ and any independence system $S$ can be approximated in polynomial time
with constant approximation ratio which is independent of $n$.

\vskip.2cm
In Section \ref{apx_sec} we prove the following result (see therein for the precise statement).

\vskip.2cm\noindent{\bf Theorem \ref{logaritmic_ratio}.}
SCO for any $c$ and any independence system $S$ can be approximated in polynomial time
with approximation ratio which is decreasing logarithmically in $n$.

\vskip.2cm
In Section \ref{small} we establish better approximation ratios
for the small values $n=2,3,4$.

\vskip.2cm
We conclude our article in Section \ref{remarks} with some final remarks.

\section{The complexity of shifted matching}\label{matchings}

In this section we discuss the computational complexity of the
shifted problem over matchings. It is convenient to introduce first
the following decision problem, of interest on its own right.

\vskip.2cm\noindent{\bf Prescribed Congestion Problem.}
Given $S\subseteq\{0,1\}^d$, $n$ and $C_1,\dots,C_d\subseteq\{0,1,\dots,n\}$,
decide if there is an $x\in S^n$ whose congestion at $i$ satisfies
$m(i,x):= \sum_{j=1}^nx_{i,j}\in C_i$ for all $i$.

\bl{congestion}
The prescribed congestion problem reduces to shifting with $c\in\{-1,0,1\}^{d\times n}$.
\el
\bpr
Define functions $f_i$ on $\{0,1,\dots,n\}$ by $f_i(j):=0$ if $j\in C_i$ and $f_i(j):=-1$
if $j\not\in C_i$. Now define $c\in\{-1,0,1\}^{d\times n}$ by $c_{i,j}:=f_i(j)-f_i(j-1)$
for $i=1,\dots,d$ and $j=1,\dots,n$.

Consider any $x\in S^n$. We then have
$$c\xx = \sum_{i=1}^d c_i\xx_i =
\sum_{i=1}^d\left(\sum_{j=1}^{m(i,x)} c_{i,j}\right)
 = \sum_{i=1}^d\left(f_i\left(\sum_{j=1}^n x_{i,j}\right)-f_i(0)\right)
 \leq -\sum_{i=1}^d f_i(0) = - |\{i\,:\,0\not\in C_i\}|$$
with equality if and only if $\sum_{j=1}^n x_{i,j}\in C_i$ for all $i$.
So the problem reduces to shifting.
\epr

The following theorem summarizes the complexity of the shifted matching problem.
\bt{matching}
Consider the shifting problem $\max\{c\xx\,:\,x\in S^n\}$ over the independence\break
system $S\subseteq\{0,1\}^d$ of matchings in a
given graph $G$ with $d$ edges and $c\in\Z^{d\times n}$. We have:
\begin{enumerate}
\item
For bipartite graphs and $c=\cc$ shifted the problem is polynomial time solvable for all $n$.
\item
For cubic graphs the problem is NP-hard already for $c=\cc$ shifted and $n=2$.
\item
For bipartite graphs and arbitrary $c$ the problem is NP-hard already for $n=2$.
\end{enumerate}
\et
Before presenting the proof of the theorem, we will show the next lemma.
\begin{lemma}\label{complexity-lem}
For any matrix $z$, let $|z|:=\sum_i\sum_j|z_{i,j}|$ and let $k\in [d]$ be an integer. Let $T\subseteq\{t\in\{0,1\}^d\,:\,\sum_{i=1}^dt_i=k\}$
and let $S:=\{s\in\{0,1\}^d\,:\,s\leq t\ \mbox{for some}\ t\in T\}$ be the independence system
generated by $T$. Then, the shifting problem $\max\{c\xx\,:\,x\in T^n\}$ over $T$ with
$c\in\Z^{d\times n}$ reduces to the shifting problem $\max\{b\xx\,:\,x\in S^n\}$
over $S$ with $b\in\Z^{d\times n}$ defined by $b_{i,j}:=c_{i,j}+2|c|+1$ for all $i,j$.
\end{lemma}
\bpr
Consider any $x,y\in T^n$ and any $z\in S^n\setminus T^n$. Then
$$b\sh x\ =\ c\sh x+(2|c|+1)|\sh x|\ =\ c\sh x+(2|c|+1)nk\ \geq\ -|c|+(2|c|+1)nk\ ,$$
$$b\sh z\ =\ c\sh z+(2|c|+1)|\sh z|\ \leq\ |c|+(2|c|+1)(nk-1)\ .$$
So $b\sh x>b\sh z$ hence an optimizer of $b\sh u$ over $S^n$ will be attained at $T^n$.
Also, $b\sh x-b\sh y=c\sh x-c\sh y$ and so an optimal solution of $b\sh u$ over $S^n$
will also be an optimal solution of $c\sh u$ over $T^n$.
\epr

We now return to prove Theorem \ref{matching}.

\bpr
Part 1 is proven in \cite{KOS}.  Note that the claim similar to Part 2 regarding perfect matchings is shown in \cite{LO}, here, for completeness, we will prove it for matchings.
Let $G$ be a cubic graph with $d$ edges. It is NP-complete
to decide if $G$ is $3$-edge-colorable \cite{Hol}. This is equivalent to deciding
if $G$ has two edge-disjoint perfect matchings. This is equivalent to the
prescribed congestion problem with $T\subset\{0,1\}^d$ the set of perfect matchings in $G$,
$n=2$ and $C_i=\{0,1\}$ for all $i$. This problem reduces by
Lemma \ref{congestion} to shifting over $T$ with shifted matrix $c=\cc$
given by $c_{i,1}=0$ and $c_{i,2}=-1$ for all $i$. By Lemma \ref{complexity-lem}, this reduces
to shifting over the independence system $S$ of all matchings in $G$ with suitable
$b$ which is shifted since $c$ is. So this problem is NP-hard and Part 2 follows.

Finally, we prove Part 3. Let $F_1,\dots,F_m\subseteq[k]$ be sets with all $|F_i|=3$.
It is NP-complete to decide if $[k]$ is partitionable by the $F_i$,
that is, if $\uplus_{i\in I}F_i=[k]$ for some $I\subseteq[m]$, see \cite{DF}.

Construct a graph $G$ with $d:=12m$ edges and $6m+2k$ vertices as follows.
For $i=1,\dots,m$ include a $6$-cycle $(u_{i,1},v_{i,1},u_{i,2},v_{i,2},u_{i,3},v_{i,3})$, which will be referred to as a {\em hexagon}.
For each $j\in[k]$ include two vertices $a_j,b_j$. For any $i\in[m]$,
 $F_i=\{r,s,t\}$ for some $1\leq r<s<t\leq k$. We introduce the six edges $\{a_r,u_{i,1}\},\{a_s,u_{i,2}\},\{a_t,u_{i,3}\}$, $\{b_r,v_{i,1}\},\{b_s,v_{i,2}\},\{b_t,v_{i,3}\}$ to the graph, these edges will be referred to as {\em non-hexagon edges}. This graph is bipartite with the $u_{i,r},b_j$ on one side and the $v_{i,r},a_j$ on the other side. We claim that the partitioning problem reduces to the prescribed congestion problem with $T\subset\{0,1\}^d$ the set of perfect matchings in $G$, $n=2$, and $C_e=\{0,1\}$ for each hexagon edge and $C_e=\{0,2\}$ for each non-hexagon edge.

Suppose $I\subseteq[m]$ gives a partition $\uplus_{i\in I}F_i=[k]$.
Construct two perfect matchings $M_1,M_2$ in $G$ as follows.
For each $i\in I$ include in both $M_1,M_2$ the six non-hexagon edges incident
on the vertices $u_{i,r},v_{i,r}$. For each $i\not\in I$ include the three hexagon
edges $\{u_{i,r},v_{i,r}\}$ in $M_1$ and the other three hexagon edges in $M_2$.
Clearly each vertex $u_{i,r},v_{i,r}$ is incident on exactly one edge in $M_1$
and one in $M_2$. Now consider $j\in[k]$. Then $j\in F_i$ for exactly one $i\in I$.
So for exactly one $1\leq r\leq 3$ we have that $a_j$ is incident in both $M_1,M_2$
exactly once on $\{a_j,u_{i,r}\}$ and $b_j$ is incident in both $M_1,M_2$ exactly
once on $\{b_j,v_{i,r}\}$. So both $M_1,M_2$ are perfect matchings in $G$,
and by the construction, the congestion on each edge $e$ is in $C_e$.

Conversely, suppose $M_1,M_2$ are two perfect matchings in $G$ with the prescribed
congestion on each edge. Construct $I\subseteq[m]$ as follows. Consider any $i\in[m]$
and the corresponding hexagon. Consider any two consecutive hexagon edges, say $\{u_{i,1},v_{i,1}\},\{v_{i,1},u_{i,2}\}$. We claim that both must have the same
congestion under $M_1 \cup M_2$. Indeed, if one has congestion $0$ and the other has congestion $1$,
then the congestion of the non-hexagon edge incident on $v_{i,1}$ must be $1$, which is impossible. So either all edges of hexagon $i$ have congestion $0$ and all non-hexagon
edges touching it have congestion $2$, in which case we include $i$ in $I$,
or all edges of hexagon $i$ have congestion $1$ and all non-hexagon edges
touching it have congestion $0$, in which case we exclude $i$ from $I$.
Next, we prove that $I$ is a solution to the partitioning problem.
For any $j\in[k]$, consider any edge $\{a_j,u_{i,r}\}$ incident
on $a_j$; then $j\in F_i$. Now either $\{a_j,u_{i,r}\}$ has congestion $0$ under $M_1,M_2$
in which case all edges of hexagon $i$ have congestion $1$ and $i\not\in I$, or
this edge has congestion $2$ in which case all edges of hexagon $i$ have
congestion $0$ and $i\in I$. Since exactly one edge incident on $a_j$ has
congestion $2$, it follows that $j$ is in exactly one $F_i$ with $i\in I$.
So $\uplus_{i\in I}F_i=[k]$ is a partitioning.

Now this prescribed congestion problem reduces by Lemma \ref{congestion} to
the shifting problem over $T$. By Lemma \ref{complexity-lem}, this reduces in turn to
shifted combinatorial optimization over the independence system $S$ of
all matchings in $G$. So this problem is NP-hard and Part 3 follows.
\epr

\section{Constant approximation of monotone shifting}
\label{constant_appx}

A matrix $x\in\{0,1\}^{d\times n}$ is {\em orthogonal} if its columns are pairwise
orthogonal, that is, have disjoint supports, which is equivalent to
$\sum_{j=1}^n x^j\in\{0,1\}^d$. In our approximation algorithms here and in
Section \ref{apx_sec} we will use the following problem over the orthogonal matrices in $S^n$.

\vskip.2cm\noindent{\bf Disjoint Union Problem (DUP).}
Given $S\subseteq\{0,1\}^d$, positive integer $n$, and $w\in\Z^d$, solve
\begin{equation}
\max\{w\sum_{j=1}^n x^j\ :\ x\in S^n\,,\ \ \sum_{j=1}^n x^j\in\{0,1\}^d\}\ .
\end{equation}

Here and in Section \ref{apx_sec} we will assume the existence of a polynomial time $\beta$-approximation algorithm
for DUP for some $\beta\leq 1$. For $\beta=1$ this means that DUP can be solved to optimality in polynomial time;
this holds for instance for the set of indicators of independent sets of a matroid, the set of indicators
of common independent sets in the intersection of two so-called {\em strongly-base-orderable matroids}, or the
set of $\{0,1\}$-valued solutions to a system $Ax\leq b$ with $A$ a totally unimodular matrix and $b$
an integer vector. However, DUP is generally NP-hard (for instance, when $S$ is the set of indicators of matchings
in a graph, already for $n=2$, see proof of Theorem \ref{matching}); so we will need the next lemma which follows
from the classical result of \cite{CFN} on the greedy algorithm for the so-called {\em maximum coverage problem}.

\bl{greedy}
Given an independence system $S\subseteq\{0,1\}^d$ presented by a linear optimization oracle, positive integer $n$, and $w\in\Z^d$, we can find in polynomial time an orthogonal matrix
$x\in S^n$ which provides a $\beta$-approximation for DUP with ratio
$\beta=1-(1-\frac{1}{n})^n \geq 1- \frac{1}{e}$.
\el
\bpr
The maximum coverage problem is to find $x\in S^n$ maximizing
$\sum\{w_i\,:\,\sum_{j=1}^n x_{i,j}\geq 1\}$.
{\em The greedy algorithm} for this problem is the following. Set $w^0:=w$ and $x^0:=0$.
For $k=1,\dots,n$ do: define $w^k\in\Z^d$ by $w^k_i:=w^{k-1}_i$ if $x^{k-1}_i=0$ and $w^k_i:=0$ if $x^{k-1}_i=1$; query the linear optimization oracle of $S$ on $w^k$ and obtain $x^k$.
Let $x:=[x^1,\dots,x^n]\in S^n$. Then the classical result of \cite{CFN} on the
greedy algorithm guarantees that $x$ has the claimed approximation ratio for the maximum coverage problem, see also \cite[page 136]{Dorit_book} for a short illuminating analysis.
Now let $y$ be the orthogonal matrix obtained from $x$ by zeroing out, in each row $i$
with $\sum_{j=1}^n x_{i,j}\geq 1$, all entries but the first $1$, so
that $\sum_{j=1}^n y_{i,j}=1$. Since $S$ is an independence system we have $y\in S^n$
and clearly for all orthogonal $z\in S^n$ we have
$$w\sum_{j=1}^n y^j=\sum\{w_i:\sum_{j=1}^n y_{i,j}\geq 1\}=\sum\{w_i:\sum_{j=1}^n x_{i,j}
\geq 1\}\geq \beta\sum\{w_i:\sum_{j=1}^n z_{i,j}\geq 1\}=\beta w\sum_{j=1}^n z^j$$
so $y$ is the desired greedily computable $\beta$-approximation for the disjoint union problem.
\epr

We proceed to show that for any shifted $c\in\Z^{d\times n}$, that is,
$c=\sh c$ having nonincreasing rows, we can use the above to approximate $\max\{c\sh x:x\in S^n\}$
to within the same approximation ratio $\beta$ in polynomial time.
Below we use some notation introduced in \cite{LO}.

Define the {\em $n$-lift} of a set of vectors $S\subseteq\{0,1\}^d$
to be the set of matrices
$$\uparrow_n\!S\ :=\ \{x\in\{0,1\}^{d\times n}\ :\ \sum_{j=1}^n x^j\in S\}\ .$$

Clearly, if $S$ is an independence system then so is $\uparrow_n\!S$.
Moreover, the following holds.

\bl{lift_lp}
A linear optimization oracle for $\uparrow_n\!S$ is efficiently realizable from one for $S$.
\el
\bpr
Let $c\in\Z^{d\times n}$ be given. For $i=1,\dots,d$ let $j(i)$ and $w_i$ be such that
$w_i=c_{i,j(i)}=\max\{c_{i,1},\dots,c_{i,n}\}$. Query the oracle of $S$ on $w=(w_1,\dots,w_d)$
and let $s\in S$ be its answer. Then $x\in \uparrow_n\!S$ defined
by $x_{i,j}:=s_j$ if $j=j(i)$ and $x_{i,j}:=0$ otherwise is optimal for $c$.
\epr

Define the {\em $n$-disjoint-union} of a set $M\subseteq\{0,1\}^U$ where $U$
is any finite ground set to be
$$\vee_n M\ :=\ \{x\in\{0,1\}^U\ :\ \exists x_1,\dots,x_n\in M\,,\ x=\sum_{k=1}^n x_k\}\ .$$

For $M\subseteq\{0,1\}^{d\times n}$ let $[M]$ be the set of
matrices equivalent to some matrix in $M$,
$$[M]\ :=\ \{x\in\{0,1\}^{d\times n}\ :\ \exists\ y\in M\,,\ x\sim y \}\ .$$

The next lemma is form \cite{LO} and its proof can be found therein.
\bl{partition}
For any set $S\subseteq\{0,1\}^d$ and any $n$ we have that
$[S^n]=\vee_n \uparrow_n\!S$ in $\{0,1\}^{d\times n}$.
\el

We can now prove the following result.
\bt{constant_ratio}
There is a polynomial time algorithm that, given a shifted $c\in\Z^{d\times n}$ and
an independence system $S\subseteq\{0,1\}^d$ presented by a linear optimization oracle,
finds $y\in S^n$ with
$$c\sh y\ \geq\ \beta\cdot\max\{c\sh z\, :\, z\in S^n\}\,,\quad\quad
\beta = 1-\left(1-\frac 1n\right)^n \geq 1-\frac 1e\ \geq\ 0.6321\ .$$
\et
\bpr
By Lemma \ref{lift_lp} we can realize a linear optimization oracle for the independence system
$M:=\uparrow_n\!S\subseteq\{0,1\}^{d\times n}$ over the ground set $U:=[d]\times[n]$.
We then consider the disjoint union problem over $M$ with $n$ and $w:=c$. By Lemma \ref{greedy}
we can find matrices $x_1,\dots,x_n\in M$ with pairwise disjoint supports,
that is, with $\sum_{k=1}^n x_k\in\{0,1\}^{d\times n}$, such that
$$c\sum_{k=1}^n x_k\ \geq\ \beta\cdot \max\{c\sum_{k=1}^n v_k\ :\ v_1,\dots,v_n\in M
\,,\ \ \sum_{k=1}^n v_k\in\{0,1\}^{d\times n}\}\ .$$
For $k=1,\dots,n$ define a vector $y^k:=\sum_{j=1}^nx_k^j$ and note that
$x_k\in \uparrow_n\!S$ implies $y^k\in S$ for all $k$. Let $x$ be the matrix
$x:=\sum_{k=1}^nx_k$ and let $y\in S^n$ be the matrix with columns $y^1,\dots,y^n$ .
Consider any feasible solution $z\in S^n$ for the shifting problem. We claim the inequality
\begin{eqnarray}\label{shifted_inequality}
c\sh y\ =\ c\sh x\ \geq\ cx\ \geq\ \beta c\sh z\ .
\end{eqnarray}
First, we have
$$\sum_{j=1}^nx^j\ =\ \sum_{j=1}^n\sum_{k=1}^nx_k^j\ =\
\sum_{k=1}^n\sum_{j=1}^nx_k^j\ =\ \sum_{k=1}^ny^k\ .$$
This implies $x\sim y$ and $\sh y=\sh x$, and so the left equality of \eqref{shifted_inequality} follows. Second, $c$ is shifted, that is, nonincreasing, and so by a trivial exchange argument, the middle inequality of \eqref{shifted_inequality} follows. Third, by Lemma \ref{partition} we have
that $z\in S^n$ implies $\sh z\in[S^n]=\vee_n \uparrow_n\!S=\vee_n M$,
and hence by the choice of $x_1,\dots,x_n$,
$$cx\ \geq\ \beta\cdot\max\{cv\,:\,v\in \vee_n M\}\ \geq\ \beta c\sh z\ ,$$
which is the right inequality of \eqref{shifted_inequality}. So inequality \eqref{shifted_inequality} holds throughout which implies that $y$ is the desired
$\beta$-approximation for the shifted combinatorial optimization problem.
\epr

\section{Logarithmic approximation of shifting}\label{apx_sec}

In this section we develop our approximation algorithm for SCO with arbitrary matrix $c$.

\subsection{The cleaning procedure}

Given a feasible solution $x$ to SCO, for $i\in [d]$, we let $m(i,x)=\sum_j x_{i,j}$ be the congestion of $i$ in $x$; we let $f_i(x)=c_i \xx_i = \sum_{j=1}^{m(i,x)} c_{i,j}$ be the {\it profit of $x$ from element $i$}, and thus the goal is to find $x\in S^n$ that maximizes $c\xx=\sum_{i\in [d]} f_i(x)$; we let $P_i(x) = \max _{\ell \leq m(i,x)} \sum_{j=1}^{\ell} c_{i,j}$ be the {\it potential profit of $x$ from element $i$}, and a maximizer in the last maximum is denoted as $p(i,x)$.  Let $\sum_{i\in [d]} P_i(x)$ be the {\it total potential profit of $x$}.  Then we have the following observation that follows by a cleaning post-processing phase:

\begin{lemma}
There exists a polynomial time algorithm that takes as an input a feasible solution $x\in  S^n$ to SCO and returns a feasible solution $x'\in S^n$ such that the objective function value of $x'$ as a solution to SCO equals the total potential profit of $x$.
\end{lemma}
\bpr
We modify $x$ as follows.  We go over all elements and for each element $i\in [d]$ we compute $m(i,x)-p(i,x)$ and if this is strictly positive, we replace exactly  $m(i,x)-p(i,x)$  $1$'s in the $i$-th row of $x$ with $0$'s.  In the resulting matrix that we denote by $x'$, the congestion of $i$ is exactly $p(i,x)$ (for every $i\in [d]$).  Thus, the objective function value of $x'$ as a solution to SCO is exactly the total potential profit of $x$.  The feasibility of $x'$ follows by our assumption that $x\in S^n$, $x'\leq x$, and $S$ is downward monotone.
\epr

The procedure described in the proof of the last lemma of transforming the solution $x$ into the solution $x'$ will be named the {\it cleaning procedure applied on $x$}.

\subsection{The approximation algorithm for SCO}
The algorithm we present to approximate SCO chooses one of $\lceil \log_2 n \rceil +1$ solutions.  That is, the following procedure is applied  for $\ell =0,1,\ldots ,\lceil \log_2 n \rceil$ with $k(\ell)= \lfloor \frac {n}{2^{\ell}} \rfloor$ for $\ell\leq \log_2 n$ and $k(\ell)=1$ for $\ell >\log_2 n$. For each such value of $\ell$, we compute its corresponding $k(\ell)$ and construct a feasible solution to SCO as described below.  Last, we choose the best solution among all the solutions we computed.

In the iteration for a fixed value of $\ell$ and its corresponding $k(\ell)$ we use the $\beta$ approximation
algorithm for approximating DUP over $S$ with the positive integer $k(\ell)$ standing for the required number of columns in the output of DUP and the weight vector $w\in\Z^d$
defined as $w_i= \max _{q \leq \min\{ 2^{\ell}, n\}} \sum_{j=1}^{q} c_{i,j}$ the maximum profit that can
be achieved from element $i$ by covering it by at most $2^{\ell}$ times.
We denote by $y(\ell)$ the matrix of dimension $d \times k(\ell)$ returned by the approximation algorithm for DUP.

If $\ell \leq \log_2 n$, then we replace each column of $y(\ell)$ by $2^{\ell}$ copies of the same column and we add another $n-k(\ell)\cdot 2^{\ell} = n - \lfloor \frac {n}{2^{\ell}} \rfloor \cdot 2^{\ell} \geq 0$ zero columns to get a feasible solution for SCO and otherwise if $\ell > \log_2 n$ we have $k(\ell)=1$ and we replace the unique column of $y(\ell)$ by $n$ copies of the same column to obtain a feasible solution for SCO.  In either case we have computed a feasible solution for SCO on which we apply the cleaning procedure to obtain our $\ell$-th candidate solution $APX_{\ell}$ whose objective function value (as a solution for SCO) is denoted by $apx_{\ell}$.

Our output is the solution $APX_{\ell}$ for the value of $\ell$ for which $apx_{\ell}$ is maximized.  We denote by $apx=\max_{\ell} apx_{\ell}$ the objective function value of the solution returned by the algorithm.

\subsection{The analysis}
Next, we turn into the analysis of this algorithm.  Let $o$ be an optimal solution for SCO.  Recall that for $i\in [d]$, $m(i,o)$ is the number of columns in $o$ that cover $i$, that is, the congestion of $i$ in $o$.  For $\ell =0,1,\ldots ,\lceil \log_2 n \rceil$ we let $$O(\ell) = \sum_{i \in [d]: 2^{\ell-1} < m(i,o) \leq 2^{\ell}} f_i(o)$$ be the total reward of $o$ from elements that are covered more than $2^{\ell-1}$ times and at most $2^{\ell}$ times,  we let
$$S(\ell)=\{ i\in [d]: 2^{\ell-1} < m(i,o) \leq 2^{\ell}\} $$ be the set of elements that contribute to the value of $O(\ell)$, and for $i\in [d]$ we let $q(i)$ be the index such that $i\in S(q(i))$.

Then, we have that the value of $o$ as a solution to SCO denoted as $OPT$ satisfies $$OPT = \sum_{\ell=0}^{\lceil \log_2 n \rceil} O(\ell) .$$
In the next lemma we prove lower bounds on $apx$ in terms of $\{ O(\ell)\}_{\ell}$.
\begin{lemma}\label{hardlemma}
For every $\ell=0,1,\ldots , \lceil \log_2 n \rceil$, we have $$apx \geq apx_{\ell} \geq \beta \cdot \sum_{i\in \bigcup_{q=0}^{\ell} S(q)} f_i(o) \cdot \left( \frac 12 \right)^{\ell-q(i)+3} = \beta \cdot \left( \sum_{q=0}^{\ell} O(q)\cdot \left( \frac 12 \right)^{\ell-q+3}  \right).$$
\end{lemma}
\bpr
Note that we have $$ \sum_{q=0}^{\ell} O(q)\cdot \left( \frac 12 \right)^{\ell-q+3} = \sum_{q=0}^{\ell} \sum_{i\in S(q)} f_i(o) \cdot  \left( \frac 12 \right)^{\ell-q+3}  = \sum_{i\in \cup_{q=0}^{\ell} S(q)} f_i(o) \cdot \left( \frac 12 \right)^{\ell-q(i)+3} ,$$ and thus the last equality in the statement of the lemma holds.

Consider the $n$ columns of $o$. We choose uniformly at random $k(\ell)$ columns of $o$ and construct a feasible
solution to SCO based on these chosen columns.  First, if the chosen columns are not pairwise disjoint,
we make them disjoint by replacing some of the $1$'s by $0$'s (we replace a minimum number of such $1$'s,
and thus do not change the set of elements that are covered by at least one chosen column).
So we obtain a collection of $k(\ell)$ pairwise disjoint columns.  Second,
 we take $\min\{ n, 2^{\ell} \}$ copies of each column in the collection of disjoint columns.  If $\ell < \log_2 n$, then we augment this collection of columns by $n-k(\ell)\cdot 2^{\ell}$ zero columns.  Then we apply the cleaning procedure on the matrix of dimensions $d \times n$ we obtained, and denote the resulting matrix by $o'$.  Observe that $o'$ is a random variable, and thus for $i\in [d]$, $f_i(o')$ is also a random variable whose expected value is denoted as $E(f_i(o'))$.

Observe that the objective function value of $o'$ as a solution to SCO is exactly the objective function value of the collection of disjoint columns as a solution to DUP.  Therefore, we know that $apx \geq apx_{\ell} \geq \beta \cdot \left( \sum_{i\in [d]} E(f_i(o')) \right)$.  Thus, it suffices to show $$\sum_{i\in [d]} E(f_i(o')) \geq \sum_{i\in \cup_{q=0}^{\ell} S(q)} f_i(o) \cdot \left( \frac 12 \right)^{\ell-q(i)+3} , $$
  and thus it suffices to show that for every $i\in \cup_{q=0}^{\ell} S(q)$, we have $$E(f_i(o')) \geq f_i(o) \cdot  \left( \frac 12 \right)^{\ell-q(i)+3} .$$

Consider such $i\in  \cup_{q=0}^{\ell} S(q)$, and let $q=q(i)$.  By the cleaning procedure, we conclude that if $i$ is covered by (at least) one chosen column, then the realization of $f_i(o')$ is at least $f_i(o)$.  Thus, it suffices to show that the probability that $i$ is covered by at least one chosen column is at least $ \left( \frac 12 \right)^{\ell-q+3}$.

For $t=1,2,\ldots k(\ell)$, we define a random variable $X_t$ that equals $1$ if the $t$-th chosen column covers the element $i$ and it equals $0$ otherwise.  Let $X=\sum_{t=1}^{k(\ell)} X_t$.  Then, we need to show that the probability $Pr(X\neq 0) \geq \left( \frac 12 \right)^{\ell-q+3}$.  To do that, we will use Chernoff bound.
First, note that $$E(X_t) = Pr(X_t=1)=1- Pr(X_t=0) \geq 1-\left( 1-\frac{2^{q-1}}{n} \right) =  \frac{2^{q-1}}{n} , $$ and thus using $k(\ell)=\max\{ 1, \lfloor \frac{n}{2^{\ell}} \rfloor \} \geq \frac{n}{2^{\ell+1}}$, we have $$E(X) \geq k(\ell) \cdot \frac{2^{q-1}}{n} \geq \frac{n}{2^{\ell+1}} \cdot \frac{2^{q-1}}{n} = \left( \frac 12 \right) ^{\ell - q+2} . $$

We use Lemma 5.27 in \cite{David+David book} to conclude that
$Pr(X=0) < e^{-\left( \frac 12 \right) ^{\ell - q+2}}.$  Thus, the probability that $i$ is covered by at least one chosen column is at least $$Pr(X\neq 0) > 1- e^{-\left( \frac 12 \right) ^{\ell - q+2}}  \geq  \left( \frac 12 \right) ^{\ell - q+3} , $$ where the last inequality holds because for every $0<z<1$ we have $e^{-z} \leq 1-z+\frac{z^2}{2} < 1-\frac{z}{2}$.
\epr
We are now in position to prove the following main result.

\bt{logaritmic_ratio}
There is a polynomial time algorithm that, given any $c\in\Z^{d\times n}$ and any
independence system $S\subseteq\{0,1\}^d$ presented by a linear optimization oracle,
finds $y\in S^n$ with
$$c\sh y\ \geq\ \frac{\beta}{4\cdot\lceil \log_2 n \rceil +8}\cdot \max\{c\sh z\, :\, z\in S^n\}\  \ ,$$
where $\beta = 1-(1-\frac 1n)^n \geq 1-\frac 1e\geq 0.6321$ is the approximation ratio of the greedy
algorithm for the disjoint union problem over $S$, and if DUP is solvable in polynomial time, then $\beta =1$.
\et
\bpr
We consider the family of constraints that we proved in Lemma \ref{hardlemma}. That is, for every $\ell$ we have $apx \geq \beta \cdot \left( \sum_{q=0}^{\ell} O(q)\cdot \left( \frac 12 \right)^{\ell-q+3}  \right).$  We multiply the constraint corresponding to $\ell=\lceil \log_2 n \rceil$ by $2$, and we sum up all the constraints.  Thus, we obtain $$\left( \lceil \log_2 n \rceil +2 \right) \cdot apx \geq \beta \cdot \sum_{\ell=0}^{\lceil \log_2 n \rceil} O(\ell) \cdot \frac 14 .$$  To see that indeed the coefficient of every $O(\ell)$ in the resulting sum is $\beta  \cdot \frac 14$ note that the coefficient is the result of a sum of geometric sequence whose ratio is $\frac 12$, the largest element in the sequence is $\frac {\beta}{8}$, and the smallest element appears twice.  These properties hold for every value of $\ell$, and every such sequence of numbers has sum of $\frac {\beta}{4}$.
\epr

\section{Improved approximations for small values of $n$}\label{small}

Here, we would like to analyze the approximation ratio of our approximation algorithm of Section \ref{apx_sec}
and variants of it for small values of $n$.
Thus, we will show that the analysis can be tightened further for these values.
We note that similar improvements are possible also for larger values of $n$.
To illustrate the improvements, we consider the cases where $n=2,3,4$.
We use the approximation algorithm for DUP and Lemma \ref{greedy} for several values of $k\leq n$.

\subsection{The case $n=2$}
Assume that $n=2$, we will apply the algorithm of Section \ref{apx_sec} for $\ell=0$  with $k(0)=2$, and for $\ell=1$ with $k(1)=1$. Using the notation of Section \ref{apx_sec}, we will prove the following lower bounds on $apx$.

\begin{lemma}
We have $apx \geq apx_1 \geq \frac 12 \cdot O(0)+O(1)$ and $apx \geq apx_0 \geq \beta O(0)$.
\end{lemma}
\bpr
For the first bound, observe that for $k=1$, we can solve DUP in polynomial time using the linear optimization oracle. We pick one column of $o$ uniformly at random.  Using the proof of Lemma \ref{hardlemma}, it suffices to show that for an element $i\in [d]$ the probability that the chosen column covers $i$ is at least $\frac 12$ if $m(i,o)=1$ and $1$ if $m(i,o)=2$.  This last claim holds because we pick each of the two columns in $o$ with probability $\frac 12$.

Similarly, if we pick the ortogonal matrix obtained from $o$ as a solution to DUP with $k=2$, then we are guaranteed to get a solution to DUP of value at least $O(0)$.  Since we can approximate DUP within an approximation ratio of $\beta$, the claim follows.
\epr

Using the two inequalities in the last lemma we get
$$ apx \cdot (1+\frac{1}{2\beta}) \geq O(0)+O(1)=OPT ,$$
and thus the approximation ratio of the algorithm is $\frac{1}{1+\frac{1}{2\beta}} = \frac{2\beta}{2\beta+1}$ where $\beta = 1-(1-\frac{1}{2})^2 = \frac{3}{4}$.  Thus, the approximation ratio of the algorithm for $n=2$ is $\frac 35$.
\begin{theorem}
The approximation ratio of the algorithm for $n=2$ is $\frac 35$ .
\end{theorem}

\subsection{The case $n=3$}
Here we apply similar arguments to the case $n=2$.  We apply our approximation algorithm for $\ell=0$ with $k(0)=3$, and for $\ell=1$ with $k(1)=1$.  For the case $\ell=1$, we note that we can solve the instance of DUP in polynomial time optimally, whereas for $\ell=0$, we approximate the instance of DUP within a factor of $\beta=1-(1-\frac 13)^3=\frac{19}{27}$. Here, we modify our notation and let $O(0)=\sum_{i\in [d]: m(i,o) = 1} f_i(o)$, and $O(1)=\sum_{i\in [d]: m(i,o)\in \{2,3\} } f_i(o)$.

 We prove the following lower bounds on $apx$.

\begin{lemma}
We have $apx \geq apx_1 \geq \frac 13 \cdot O(0)+ \frac{2}{3} \cdot O(1)$ and $apx \geq apx_0 \geq \beta O(0)$.
\end{lemma}
\bpr
For the first bound, we pick one column of $o$ uniformly at random.  Using the proof of Lemma \ref{hardlemma}, it suffices to show that for an element $i\in [d]$ the probability that the chosen column covers $i$ is  $\frac 13$ if $m(i,o)= 1$ and at least $\frac 23$ if $m(i,o)\in \{ 2,3\}$.  This last claim holds because we pick each of the three columns in $o$ with probability $\frac 13$.

For the second bound, if we pick the ortogonal matrix obtained from $o$ as a solution to DUP with $k=3$, then we are guaranteed to get a solution to DUP of value at least $O(0)$.  Since we can approximate DUP within an approximation ratio of $\beta$, the claim follows.
\epr

By multiplying the inequality $apx \geq \frac 13 \cdot O(0)+ \frac{2}{3} \cdot O(1)$ by $\frac{3}{2}$ we get $\frac 32 \cdot apx \geq \frac 12 \cdot O(0) + O(1)$.  Then by adding $\frac 1{2\beta}$ times the inequality  $apx \geq \beta O(0)$ ,we get $$\left( \frac 32 + \frac 1{2\beta} \right) \cdot apx \geq O(0)+O(1)=OPT, $$
and thus the approximation ratio of the algorithm for $n=3$ is $\frac 1{\left( \frac 32 + \frac 1{2\beta} \right)} = \frac{2\beta}{3\beta+1} = \frac{19}{42}$.
\begin{theorem}
The approximation ratio of the algorithm for $n=3$ is $\frac{19}{42}$.
\end{theorem}

\subsection{The case $n=4$}
We apply our approximation algorithm for $\ell=0$ with $k(0)=4$, for $\ell=1$ with $k(1)=2$, and for $\ell=2$ with $k(2)=1$.  For the case $\ell=2$, we note that we can solve the instance of DUP in polynomial time optimally, whereas for $\ell=1$, we approximate the instance of DUP within a factor of $1-(1-\frac 12)^2=\frac{3}{4}$, and for $\ell=0$, we approximate the instance of DUP within a factor of $1-(1-\frac 14)^4=\frac{175}{256}$. Here, we use our original notation and let $O(0)=\sum_{i\in [d]: m(i,o) = 1} f_i(o)$,  $O(1)=\sum_{i\in [d]: m(i,o)=2} f_i(o)$, and $O(2)=\sum_{i\in [d]: m(i,o)\in \{ 3,4\} } f_i(o)$.

We prove the following lower bounds on $apx$.

\begin{lemma}
We have   $apx \geq apx_2 \geq \frac 14 \cdot O(0)+ \frac{1}{2} \cdot O(1)+ \frac{3}{4} \cdot O(2)$,   $apx \geq apx_1 \geq \frac{3}{4} \cdot \left(  \frac 12 \cdot O(0)+ \frac{5}{6} \cdot O(1) \right) = \frac 38 \cdot O(0) + \frac 58 \cdot O(1)$, and $apx \geq apx_0 \geq \frac{175}{256} \cdot O(0)$.
\end{lemma}
\bpr
For the first bound, we pick one column of $o$ uniformly at random.  Using the proof of Lemma \ref{hardlemma}, it suffices to show that for an element $i\in [d]$ the probability that the chosen column covers $i$ is  $\frac 14$ if $m(i,o)= 1$, it is $\frac 12$ if $m(i,o)=2$, and it is at least $\frac 34$ if $m(i,o)\in \{ 3,4\}$.  This last claim holds because we pick each of the four columns in $o$ with probability $\frac 14$.

For the second bound, recall that we approximate the corresponding DUP instance within a factor of $\frac 34$.   We pick a pair of distinct columns of $o$ uniformly at random.  Using the proof of Lemma \ref{hardlemma}, it suffices to show that for an element $i\in [d]$ the probability that at least one chosen column covers $i$ is  $\frac 12$ if $m(i,o)= 1$ and it is $\frac {5}{6}$ if $m(i,o)=2$.  This last claim holds because we choose each pair of columns with equal probability (that is, $\frac 16$) , and there is only one pair of columns that does not cover an element $i\in [d]$ with $m(i,o)=2$.

For the last bound, if we pick the ortogonal matrix obtained from $o$ as a solution to DUP with $k=4$, then we are guaranteed to get a solution to DUP of value at least $O(0)$.  Since we can approximate DUP within an approximation ratio of $\frac{175}{256}$, the claim follows.
\epr

We multiply the constraint  $apx \geq \frac 14 \cdot O(0)+ \frac{1}{2} \cdot O(1)+ \frac{3}{4} \cdot O(2)$  by $\frac 43$, we multiply the constraint  $apx \geq  \frac 38 \cdot O(0) + \frac 58 \cdot O(1)$ by $\frac{8}{15}$, and we multiply the constraint $apx \geq  \frac{175}{256} \cdot O(0)$ by $\frac{1792}{2625}$.  Last we sum up the resulting three constraints to get that $apx \cdot \left( \frac 43 + \frac {8}{15} + \frac {1792}{2625} \right) \geq O(0)+O(1)+O(2) = OPT$.  Thus $apx \cdot \frac{3500+1400+1792}{2625} = apx \cdot \frac{6692}{2625} \geq OPT$, and we conclude the following theorem.

\begin{theorem}
The approximation ratio of the algorithm for $n=4$ is $\frac{2625}{6692} \sim 0.392259$.
\end{theorem}

\section{Remarks}\label{remarks}

First, we point out the following possible generalization of shifted optimization,
which in the congestion games context of \cite{Ros} corresponds to an individual set
of strategies for each player. We have a matrix $c\in\Z^{d\times n}$ as before,
but are now given $n$ sets $S_1,\dots,S_n\subseteq\{0,1\}^d$ rather than one. The generalized shifted optimization problem is then $\max\{c\xx\,:\,x^k\in S_k\}$. However, we note that this
problem becomes hard quickly. Let $n=3$ and let $S_1,S_2,S_3$ be matroids of common rank $r$.
Let $c_{i,1}:=c_{i,2}:=0$ and $c_{i,3}:=1$ for all $i$. Then the optimal value of
this generalized shifted problem is $r$ if and only if the three matroids have a common basis.
It is easy to construct three matroids from a given digraph that have a common basis
if and only if the digraph has an $s-t$ Hamiltonian path, which is hard to decide,
see \cite[Chapter 2]{Onn}.

Another possible generalization of the problem is the following.
We have one nonempty $S\subseteq\{0,1\}^d$ as before, but are now given a function $f:\{0,1,\dots,n\}^d\rightarrow\Z$. The objective value of $x\in S^n$ is $f$ evaluated
at the vector whose $i$-th component is the congestion of $x$ at $i$. So the problem is
$\max\{f(\sum_{j=1}^nx^j)\,:\,x\in S^n\}$. Assuming $S\neq\emptyset$ we have the following.
\bp{convex}
If $f$ is convex then there is an optimal solution $x$ having identical columns, and so the
generalized shifted problem reduces to solving $\max\{f(ns)\,:\,s\in S\}$ over the set $S$.
\ep
\bpr
Assume $f$ is convex and let $\hat x$ be an optimal solution to the generalized shifted problem. Let $s^1,\dots,s^m\in S$ be the distinct columns of $\hat x$. Define a $d\times m$ matrix $M:=[s^1,\dots,s^m]$ and a function $g:\Z^m\rightarrow\Z$ by $g(y):=f(My)$, which is convex since $f$ is. Consider the integer simplex $Y:=\{y\in\Z_+^m\,:\,y_1+\cdots+y_m=n\}$ and
the auxiliary problem $\max\{g(y)\,:\,y\in Y\}$. For each $y\in Y$ let $x(y):=[s^1,\dots,s^1,\dots,s^m,\dots,s^m]\in S^n$ be the $d\times n$ matrix
consisting of $y_i$ copies of $s^i$. Permuting the columns of $\hat x$
we may assume that $\hat x=x(\hat y)$ for some $\hat y\in Y$.

Now, $g$ is convex, so the auxiliary problem has an optimal solution which is a vertex
of $Y$, namely, a multiple $\tilde y=n{\bf 1}_i$ of a unit vector in $\Z^m$. It then follows that $\tilde x:=x(\tilde y)=[s^i,\dots,s^i]$ is the desired optimal solution for the shifted problem, proving the proposition, since
$$f(\sum_{j=1}^n\tilde x^j)\ =\ f(M\tilde y)\ =\ g(\tilde y)
\ \geq\ g(\hat y)\ =\ f(M\hat y)\ =\ f(\sum_{j=1}^n\hat x^j)\ .\mepr$$

Finally we note that our work in this article is part of a more general line of study
of the possibility of amplifying linear to nonlinear optimization
over independence systems. In \cite{LOW} this line
was also considered, for a certain family of nonlinear multicriteria problems.
It was shown therein that under some assumptions, a so termed {\em $r$-best}
solution, having the property that at most $r$ better objective values can be
attained by other solutions, can be obtained in polynomial time. It was also shown
therein that under milder assumptions, exponentially many queries to the
linear optimization oracle presenting the independence system may be required.
See \cite[Chapter 6]{Onn} and references therein for more details on this line of research.

\section*{Acknowledgment}

M. Kouteck\'y was partially supported by the project 17-09142S of the
Czech Science Foundation. A. Levin was partially supported by a Grant from
GIF - the German-Israeli Foundation for Scientific Research and Development.
S.M. Meesum was partially supported by a grant at the Technion. 
S. Onn was partially supported by the Dresner Chair at the Technion.

\end{document}